\documentclass[11.5pt,fleqn]{article}

\usepackage[english]{babel}
\usepackage{amsmath, amssymb, multicol, amsfonts, amsthm, amscd, epsfig, enumerate,xcolor}
\usepackage[latin1]{inputenc}
\usepackage[T1]{fontenc}
\usepackage{graphicx}
\usepackage{geometry}
\usepackage{float}

\newtheorem{theorem}{Theorem}

\newtheorem{corollary}[theorem]{Corollary}

\newtheorem{proposition}[theorem]{Proposition}
\newtheorem{remark}[theorem]{Remark}

\begin{document}
	
\author{ \emph{\textbf{Ramponi A.}} \\
Department of Economics and Finance, University of Rome Tor Vergata\\
Via Columbia 2, \ 00133 Rome, \ Italy\\
\texttt{alessandro.ramponi@uniroma2.it} \\ \\
\emph{\textbf{Tessitore M.E.}}\\
Department of Economics and Finance, University of Rome Tor Vergata \\
Via Columbia 2, \ 00133 Rome, \ Italy\\
\texttt{tessitore@economia.uniroma2.it} }

\title{The economic cost of social distancing during a pandemic: an optimal control approach in the SVIR model}
\date{\today}
\maketitle

\begin{abstract}
We devise a theoretical model for the optimal dynamical control of an infectious disease whose diffusion is described by the SVIR compartmental model. The control is realized through implementing social rules to reduce the disease's spread, which often implies substantial economic and social costs. We model this trade-off by introducing a functional depending on three terms: a social cost function, the cost supported by the healthcare system for the infected population, and the cost of the vaccination campaign. Using the Pontryagin's Maximum Principle, we give conditions for the existence of the optimal policy, which we characterize explicitly in three instances of the social cost function, the linear, quadratic and exponential models, respectively. Finally, we present a set of results on the numerical solution of the optimally controlled system by using Italian data from the recent Covid--19 pandemic for the model calibration.
 
\end{abstract}

\thispagestyle{empty}


\smallskip

\noindent
\textbf{AMS2020}: 49K15, 49M05, 34H05, 92B05.

\smallskip


\section{Introduction}

For almost three years, around the world,  governments have been trying to figure out the best policy to manage the pandemic caused by Covid--19. This sudden global emergency has highlighted the need to systematically address the problem of managing epidemics in a closely interconnected society. The containment measures which have been considered, from the mildest ones, such as the use of masks, to the more limiting ones, such as periods of home confinement, come naturally at the expense of losing the benefits of contact, and they may induce a persistent economic depression due to many aspects. From one side, the impossibility of carrying out all or part of the usual work activity impacts the production system as a whole. On the other hand, people become afraid of leaving their home, thus forgoing all those activities based on "social contact", such as purchasing goods, personal healthcare (e.g., preventive health checks), travelling, and so on. Although in almost all countries, to overcome some of these limitations, there has been increasing use of web applications for smart working, online shopping, and, to some extent, some kind of social activities, these measures entail some costs to society besides the ones strictly due to the disease itself, such as costs of treatments for infected people or for implementing a vaccinations campaign. The development of preventive or control interventions for the disease is therefore of central importance, as is the cost-effectiveness estimation of such measures. The definition of appropriate optimal control strategies thus becomes essential to study the planner trade-off between the direct cost due to the spread of the disease and containing the disease through social measures, which certainly implies an economic effort for the society as a whole.

In this paper, we consider the interaction between social distance interventions and the spread of a disease by using a control theoretic approach. Optimal control theory is undoubtedly the key tool to attain the trade-off between the fight against disease and the costs of social limitations imposed by a planner.

The starting point in most of the literature is the description of the spread of an infectious diseases by means of the SIR model, or its several modification (see e.g. Brauer,and Castillo-Chavez in \cite{BC}), introduced by Kenmark and Mckendrick in \cite{KM}, where each letter represents a compartment in which any individual of a population can be set: Susceptible (S), Infected (I), and Recovered/Removed (R). The literature on infectious disease analyzed via optimal control (see e.g. Lenhart and Workman in \cite{LW}) is rapidly increasing. Behncke \cite{Ben00} is one of the first attempt to use systematically a control approach in the framework of the epidemiological models. In the past decades the research was focused on measures based on selective isolation and immunization. Abakuks, in \cite{A}, assuming that an infected population can be instantaneously isolated,  studied how to optimally separate it, while Hethcote and Waltman proposed optimal vaccination strategies in \cite{HW73}. More recently, Ledzewicz and Sch\"attler, \cite{LS}, were dealing with an optimal control problem using a model with vaccines and treatments on a growing population, while Federico et al.,in \cite{FFT22}, studied an optimal vaccination strategy in a SIRS compartmental model, using a dynamic programming approach. In \cite{GS} Gaff and  Schaefer considered SIR/SEIR/SIRS models where the controls are again on the vaccination rate and the cure given to the infected persons, who could also be quarantined. In \cite{B}, Bolzoni et al. analyzed time-optimal control problem for the use of vaccination, isolation, and culling in the linear case. In Miclo et al. \cite{MSW20}, the authors consider a deterministic SIR model in which the social planner controls the transmission rate in order to lower its natural level so as not to overburden the health care system. Federico and Ferrari \cite{FF}, dealt with the issue of a policymaker aiming to tame an epidemic's spread while minimizing its associated social costs in a stochastic extension of the SIR model. 
Kruse and Strack, in \cite{KS}, extend the SIR model with a parameter controlled by the planner, which affects the rate at which the diseases are transmitted, capturing political measures such as social distancing and lockdown of institutions and businesses. While these measures reduce the spread of the disease, they often lead to economic and social costs. They model this trade-off by considering convex costs in the number of infected and the reduction in transmission rate. The control through lockdown policies, which affect the rate of diffusion of the disease in a SIRD model, is studied in Calvia et al. \cite{CGLZ22} using a dynamic programming approach. 

Undoubtedly, one of the most widely used preventive interventions is vaccination. Nowadays, there is extensive literature on vaccination models; see, for instance, the book by Brauer and Castillo--Chavez \cite{BC}. In order to include a vaccination strategy explicitly into the dynamical description of the disease, we rely on the model proposed by Liu et al. \cite{XTI}, denoted as SVIR. Indeed, they consider vaccination in a basic SIR model by introducing a new compartment $V$ where the vaccinees will belong before reaching immunity and, therefore, entering the compartment $R$ of recovered individuals. 

The application of an optimal control approach to a SVIR dynamical model is less considered in the literature. In \cite{KPS}, Kumar and Srivastava propose and analyze a control problem in this framework by using vaccination and treatment as control policies, and a cost functional linear in the state variables, quadratic in the treatment and quartic in the vaccination policies, respectively. Witbooi et al. in \cite{WMV}, considered both a deterministic and stochastic optimal problem for the SVIR model, assuming the vaccination rate as control, and an additive cost functional. 

In this paper, we assume a deterministic SVIR dynamical model to describe the spread of an infectious disease that a social planner can control through a set of mitigation measures whose aim is to lower the rate of contagion in the population. The challenge is to find the optimal response balancing restrictions that will minimize the incidence of the disease, keeping in mind the economic cost of such limitations and having at disposal an immunization instrument. We, therefore, introduce an explicit cost function to take into account the impact of such measures other than the cost of vaccination and the cost due to the infected population, then proving the existence of a solution for the optimal control problem by using the Pontryagin's Maximum Principle (see, e.g., \cite{LW}). Hence, by specifying the functional form of the social cost function, the linear, quadratic and exponential instances,  we are able to explicitly derive the optimal control strategy function. The obtained results are analyzed through a set of numerical experiments based on the implementation of the Forward-Backward Sweep algorithm (\cite{LW}), by calibrating some of the model parameters using the Italian dataset of the recent Covid--19 pandemic.   

As expected, by setting costs for treatment and vaccination, as the social cost increases, the optimal strategy is pushed toward the no-control strategy, while conversely, if the social cost is negligible, the optimal strategy becomes the maximum-control strategy. By fixing a comparable value for the social cost, the optimal control, in the quadratic and exponential case, is typically at its maximum level in the first period of the pandemic, then it decreases to a lower value when the level of infected compartment becomes small. Interestingly, in the linear case, which implies a bang-bang type of optimal control, after a period of maximum control, which corresponds to a decreasing behavior of the infected population, the optimal control goes to the minimum level involving a temporary increase in the number of infected, which subsequently tend to decrease again. Therefore, the I compartment dynamic shows a "waves" behaviour.
We finally analyzed the optimal strategy in an example having an endemic equilibrium.

Our paper is structured as follows:  first of all, in Section 2 we recall the basic SVIR model together with its main properties, then in Section 3 we formulate the deterministic optimal control problem, proving the existence of a solution, and characterizing explicitly the optimal control for several instances of the cost functional. Finally, by using the Forward-Backward sweep algorithm we numerically solve the problem, and we illustrate the results obtained in several examples.

\section{The Basic SVIR Model} \label{basic_svir}

The SVIR model was introduced by Liu et al. in \cite{XTI} to modify the well-known SIR model in order to include a vaccination program (continuous or impulsive) in the considered population. The four groups are, therefore, the Susceptibles $S$, the Infected $I$, the Recovered $R$, and the Vaccinees $V$, representing those having begun the vaccination process, where $S$, $V$, $R$ and $I$ denote the fractions of the total population belonging to each group, respectively.

\begin{figure}[h]
\begin{center}
 	\includegraphics[width=5cm, height=3cm]{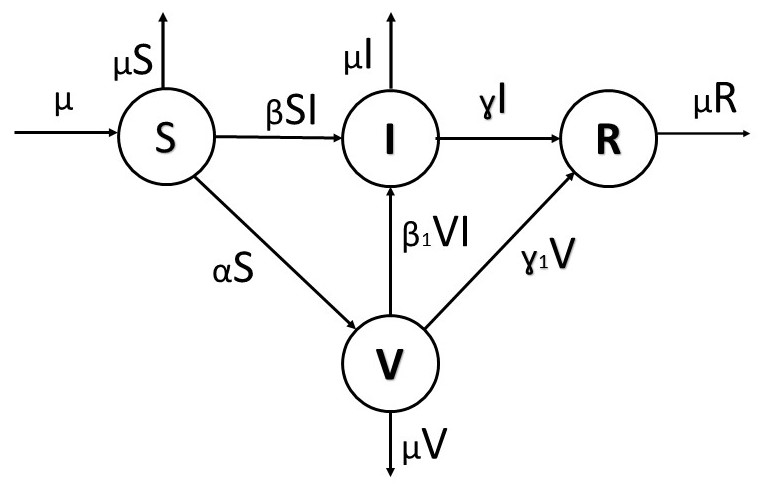}
\end{center}
\caption{The basic SVIR model graph.}
\label{svir_grafo1}
\end{figure}

Let $\beta $ represent the transmission rate of disease when the susceptible individuals get in contact with the infected ones, and let $\gamma $ be the recovery rate of the infected individuals. It is assumed that vaccinated individuals gain immunity against the disease at a rate $\gamma _1$  and that even the vaccinees have the chance to be infected at a rate $\beta _1$, which can be taken smaller than $\beta $ since after the vaccination process some immunity is acquired. Parameter $\alpha $ represents the rate at which the susceptible persons are moving in the vaccination program, and $\mu $ is the birth-death rate. Figure \ref{svir_grafo1} shows how the population is moving among  the four compartments $S,V,I,R$. 

The framework for the continuous vaccination process can be described through the following system of first order differential equations:
\begin{equation}
\label{svirconst}
\left\{
\begin{array}{ll}
\displaystyle{\frac{dS}{dt}(t)} =-\beta S(t)I(t) -\alpha S(t)+ \mu - \mu S(t) &S(0)=S_0
\\\\
\displaystyle{\frac{dV}{dt}}(t) =\alpha S(t) - \beta_1  V(t)I(t)-\gamma _1V(t)- \mu V(t) &V(0)=V_0
\\\\
\displaystyle{\frac{dI}{dt}}(t) =\beta S(t)I(t) + \beta_1 V(t)I(t)-\gamma I(t) - \mu I(t) &I(0)=I_0
\\\\
\displaystyle{\frac{dR}{dt}}(t)=\gamma _1V(t)+\gamma I(t)- \mu R(t) &R(0)=R_0
\end{array}
\right.
\end{equation} 
where the parameters $\beta, \beta_1, \gamma, \gamma_1, \mu \in \mathbb{R}^+$  and $\alpha \geq 0$. 
Moreover,  we assume that the initial data   $S_0, V_0, I_0, R_0 \in \mathbb{R}^+$, and $S_0+V_0+I_0+R_0 = 1$. The above assumptions are stated since the model (\ref{svirconst}) represents human populations, and it can be shown that the solutions of the system are non--negative given non--negative initial values, see \cite{XTI}. In particular, it is worth noticing that by defining $N(t) = S(t)+V(t)+I(t)+R(t)$, we immediately have from (\ref{svirconst}) that $\frac {dN}{dt}(t) = 0$: hence $N(t)=N_0 \equiv 1$, for all $t\geq 0$.

Since the last equation in the system (\ref{svirconst}) is a linear combination of the previous ones, it is enough to study the properties of the system using only the three state variables $S$, $V$, and $I$. In \cite{XTI} it is shown that the model (\ref{svirconst}) has a disease free equilibrium (that is an equilibrium $(S^*, V^*, I^*)$ for which $I^* \equiv 0$) 
\begin{equation}
E_0=\left(\frac{\mu}{\mu+\alpha}, \frac{\alpha \mu}{(\mu+\gamma_1)(\mu+\alpha)}, 0 \right)
\end{equation}
and an endemic equilibrium 
\begin{equation} 
E_+ = \left(\frac{\mu}{\mu+\alpha+\beta I_+}, \frac{\alpha \mu}{(\mu+\alpha+\beta I_+)(\mu+\gamma_1+\beta_1 I_+)}, I_+\right),
\end{equation}
where $I_+$ is the positive root of quadratic equation, whose coefficients depend on the parameters of the model and on the basic reproduction number, given by
\begin{equation}\label{R0}
R_0^C = \frac{\mu \beta}{(\mu+\alpha)(\mu+\gamma)}+\frac{\alpha \mu \beta_1 }{(\mu+\gamma_1)(\mu+\alpha)(\mu + \gamma)}.
\end{equation}

The main properties of the dynamical system (\ref{svirconst}) are summarized in the two following Theorems:

\begin{theorem}
The disease free equilibrium $E_0$, which always exists, is locally asymptotically stable if $R^C_0 <1$ and is unstable if $R^C_0>1$. System (\ref{svirconst}) has a unique positive equilibrium $E_+$ if and only if $R^C_0 >1$ and it is locally asymptotically stable when it exists.
\end{theorem}

\begin{theorem}
If $R^C_0 \leq 1$, then the disease free equilibrium $E_0$ is globally asymptotically stable. And if $R^C_0 > 1$, the endemic equilibrium $E_+$ is globally asymptotically stable.
\end{theorem}


\section{The controlled SVIR Model} 
\label{sect_control}

In this Section, we introduce the controlled SVIR model, and we analyze a deterministic optimal control problem associated with it.

If an infectious disease occurs, we consider a control variable $u(\cdot )$, which is meant to govern the social restrictions imposed by the ruler on a population until a specific time $T$, which is the final time of government restrictions. Even if $T$ itself could be a decision variable, we prefer in this paper to consider it as fixed in advances. 

The control variable $u$ belongs to the admissible set $U$ defined as 
$$U = \{u(t):\mathbb{R}^+\to [0,1]| \; 0 \leq u(t) \leq \bar u,\; t \in[0,T], \mbox{ Lebesgue measurable}\}, \mbox{ where }\bar u\in [0,1].$$

The control variable $u$ allows to \textit{adjust} the rate of transmission of the disease, which we model as a decreasing linear function $\beta(\cdot)$. 
We want to design the situation where in the absence of control ($u \approx 0$), the infectivity rate $\beta$ is high, while for increasing controls, this rate decreases. The function $\beta $ captures both the infectivity of the disease and the restrictions ruler imposes to govern the speed at which the infection spreads. 
Thus we get the following controlled SVIR model (see Fig. \ref{contr_svir_grafo}):
\begin{figure}[h]
    \centering
    \includegraphics[width=5cm, height=3cm]{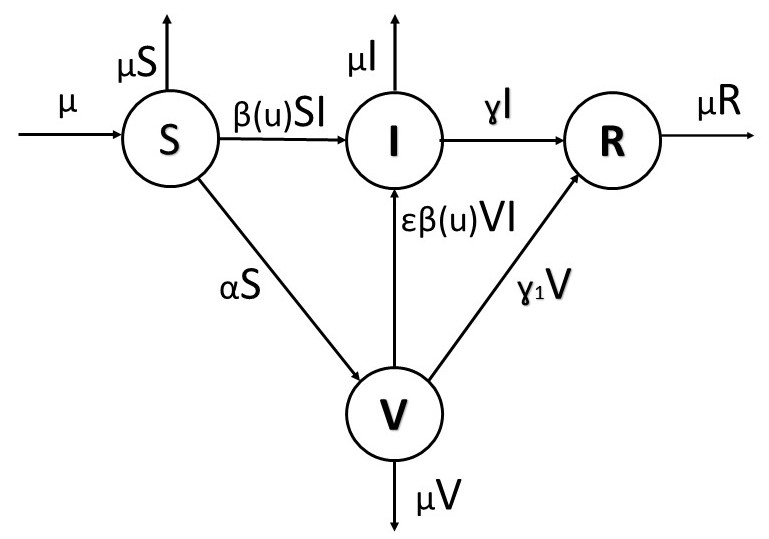}
    \caption{The controlled SVIR model graph.}
    \label{contr_svir_grafo}
\end{figure}

\begin{equation}
\label{sistema}
\left\{
\begin{array}{ll}
\displaystyle{\frac{dS}{dt}(t)} =-\beta (u(t))S(t)I(t) -\alpha S(t)+ \mu - \mu S(t) &S(0)=S_0
\\\\
\displaystyle{\frac{dV}{dt}}(t) =\alpha S(t)-\varepsilon \beta (u(t))V(t)I(t)-\gamma _1V(t)- \mu V(t) &V(0)=V_0
\\\\
\displaystyle{\frac{dI}{dt}}(t) =\beta (u(t))S(t)I(t) +\varepsilon \beta (u(t))V(t)I(t)-\gamma I(t) - \mu I(t) &I(0)=I_0
\\\\
\displaystyle{\frac{dR}{dt}}(t)=\gamma _1V(t)+\gamma I(t)- \mu R(t) &R(0)=R_0
\end{array}
\right.
\end{equation}
where $\varepsilon$ quantifies the vaccine efficacy (if $\varepsilon \equiv 0$ no vaccinated gets infected). Furthermore, we assume that $S_0, V_0, I_0, R_0 \in \mathbb{R}^+$.

Since the last equation of (\ref{sistema}) is a linear combination of the other equations, and the recovered people will not enter in the specification of the costs of the disease, it is enough to consider the following system

\begin{equation}
\label{sistemar}
\left\{
\begin{array}{ll}
\displaystyle{\frac{dS}{dt}}(t) =-\beta (u(t))S(t)I(t) -\alpha S(t)+ \mu - \mu S(t) &S(0)=S_0
\\\\
\displaystyle{\frac{dV}{dt}}(t) =\alpha S(t)-\varepsilon \beta (u(t))V(t)I(t)-\gamma _1V(t)- \mu V(t) &V(0)=V_0
\\\\
\displaystyle{\frac{dI}{dt}}(t) =\beta (u(t))S(t)I(t) +\varepsilon \beta u(t)V(t)I(t)-\gamma I(t) - \mu I(t) &I(0)=I_0
\end{array}
\right.
\end{equation}

Now we can formulate the optimal control problem. To this end, we introduce the following functional in order to minimize the cost of the infected population $I$ and the cost of the vaccination, which is proportional to a fraction $\alpha$ of the susceptible individuals in $S$. We suppose that these costs are due to hospitalization expenses for patients requiring inpatient care, with or without ICU (Intensive Care Unit), and to the arrangement of the vaccination program supply chain (e.g., the setting up and the management of a vaccination hub, of the related medical staff, and so on). Moreover, we assume that the cost of social restrictions is a function $c$ of the control variable $u$ such that  $c$ is a strictly increasing, strictly convex function of the restriction policy $u$, and that $c(0)=0$. This means that, in the absence of control, the total costs of the disease diffusion are due to the infected individuals and the vaccination strategy. In this way, and by assuming an additive structure for the cost functional, we disentangle the costs entirely due to the disease from those due to the ``restrictions'' imposed on the whole society. Parameters $c_1, c_2 \in \mathbb{R}^+$ represent the cost of being infected and the cost of the vaccination campaign respectively.

Hence the objective function is given by 
\begin{equation}
\label{funzionale}
J(u) =\int_{0}^{T} [c (u(t))+c_1I(t)+c_2\alpha S(t)]dt.
\end{equation}
Our goal is to derive the optimal strategy $u^*\in U$ and the associated state variables $S, V$ and $R$ to minimize (\ref{funzionale}). To prove the existence of such a strategy $u^*$ we refer to \cite{FR}, \cite {LW}
and  \cite{KPS}.

\begin{theorem}
Let $\beta (\cdot)$ be a linear decreasing function and let $c(\cdot)$ be a twice continuous differentiable function, such that $c'>0$, $c''>0$ and $c(0)=0$. 

Then an optimal solution $ u^*$  for problem (\ref{sistemar})--(\ref{funzionale}) exists, i.e. there exists an optimal control $u^*\in U$ such that $J(u^*)=\min J(u)$.
\end{theorem}

In order to solve the above optimal control problem, we refer to the well-establish control theory, see for instance \cite{FR} or \cite{LW}.
It is introduced the Hamiltonian function $H$ and the Lagrange multipliers $\lambda _1(\cdot )$, $\lambda _2(\cdot )$ and $\lambda _3(\cdot )$ denoted also as co--state variables. From now on, even if the state variables $S, V, I$, the control variable $u$ and the co-state variables $\lambda _1$, $\lambda _2$ and $\lambda _3$ are functions of time, we omit this dependence except where it is explicitly required. 

The Hamiltonian function of the optimal control problem (\ref{sistemar})--(\ref{funzionale}) is defined as follows

\begin{equation}
\label{hamiltoniana}
\begin{array}{l}
H(t,S,V,I,u) 
=c(u)+c_1I+c_2\alpha S+
\lambda _1[-\beta (u)SI -\alpha S+ \mu - \mu S ]+
\\\\
+\lambda _2[\alpha S-\varepsilon \beta (u)VI-\gamma _1V- \mu V]
+\lambda _3[\beta (u)SI +\varepsilon \beta (u)VI-\gamma I - \mu I].
\end{array}
\end{equation}

\begin{theorem} \label{optimal_system}
An optimal solution $(S^*,V^* , I^*, u^*)$ for problem (\ref{sistemar})--(\ref{funzionale})  satisfies the following system of differential equations
\begin{equation}
\label{costato}
\left\{
\begin{array}{l}
\lambda _1' =[\beta (u)I +\alpha + \mu]\lambda _1 - 
\alpha \lambda _2-\beta (u)I\lambda _3-c_2\alpha 
\\\\
\lambda _2' =[\varepsilon \beta (u)I +\gamma _1+ \mu]\lambda _2 
-\varepsilon \beta (u)I\lambda _3
\\\\
\lambda _3' =\beta (u)S \lambda _1+\varepsilon \beta (u)V\lambda _2-
[\beta (u)S+\varepsilon \beta (u)V-\gamma -\mu]\lambda _3-c_1
\end{array}
\right.
\end{equation}
with the transversality conditions on the co--states $\lambda _1$, $\lambda _2$ and $\lambda _3$ given by

$$\lambda _1(T)=0\qquad \lambda _2(T)=0\qquad \mbox{and}\qquad \lambda _3(T)=0$$
The optimal restriction policy $u^*$ is such that
\begin{equation}
\label{u*}
\displaystyle{
u^*(t)\in \mbox{argmin}_{u\in [0,1]} H(t,S,V,I,u) }.
\end{equation}

\end{theorem}
\proof 
Let  $(S^*,V^* , I^*, u^*)$ be an optimal solution for problem (\ref{sistemar})--(\ref{funzionale}).  By Pontryagin's Maximum Principle the costate  variables $\lambda _1$, $\lambda _2$ and $\lambda _3$   satisfy system (\ref{costato} ) whose equations are obtained evaluating the partial derivatives of the Hamiltonian function $H$ (\ref{hamiltoniana}), with respect to the state variables $S,V,I$  
\begin{equation}
\left\{
\begin{array}{l}

\lambda _1' =-\displaystyle{\frac{\partial H}{\partial S}}
\\\\
\lambda _2' =-\displaystyle{\frac{\partial H}{\partial V}}
\\\\
\lambda _3' =-\displaystyle{\frac{\partial H}{\partial I}}
\end{array}
\right.
\end{equation}
 with the transversality conditions 
 $\lambda _1(T)=\lambda _2(T)=\lambda _3(T)=0.$
The Hamiltonian function $H$, defined in (\ref{hamiltoniana}) is strictly convex with respect to the control variable $u$, hence the existence of a unique minimum follows, see \cite{WMV}, hence 
$$
u^*(t)\in \mbox{argmin}_{u\in [0,1]} H(t,S,V,I,u).
$$
\qed

\begin{remark}
Similar results are easily generalized using a convex function for the cost of the infected population. This more general assumption models the non-linear impact of disease spread on the healthcare system, resulting in hospital services becoming overwhelmed.
\end{remark}

\bigskip

We now further specialize the result obtained by explicitly specifying the functional form of the transmission rate $\beta(u)$ and cost function $c(u)$. As a basic model we consider the following linear model: 
\begin{equation} \label{lin_contr}
\beta(u) = \beta_0 (1- u), \ \ \ 0 \leq u \leq 1,
\end{equation}
where $\beta_0>0$ is the specific transmission rate of the disease. In this case, we model the situation when the maximum control (i.e. $u \equiv 1$) completely ``freeze'' the disease diffusion.

In our practical application we consider the following functions:
\begin{enumerate}
\item $c_Q(u) = b u^2$, $ b >0$;
\item $c_{exp}(u) = e^{k u}-1$, $k >0$;
\item $c_{lin}(u)=au$, $a>0$. 
\end{enumerate}

A complete characterization of the optimal controls is proved in the following.

\begin{corollary}
Let $\beta (u) = \beta_0 (1- u)$ and $c_Q(u)=b u^2$. Then the optimal control strategy $u_{Q}^*$ for  problem (\ref{sistemar})--(\ref{funzionale}) is given by
\begin{equation}
\label{optcontrQdr}
\displaystyle{u_{Q}^*(t)=min\left \{max \left [0,\frac{\beta_0I(t)[S(t)(\lambda _3(t)-\lambda _1(t))+
\varepsilon V(t)(\lambda _3(t)-\lambda _2(t))]}{2b}\right ],\overline{u}\right \}}
\end{equation}
\end{corollary}
\proof
In this case the Hamiltonian function $H$ is defined as
\begin{equation}
\label{hamiltoniana2}
\begin{array}{l}
H(t,S,V,I,u) 
=bu^2+c_1I+c_2\alpha S+
\lambda _1[-\beta_0(1- u) S I -\alpha S+ \mu - \mu S ]+
\\\\
+\lambda _2[\alpha S-\varepsilon \beta_0 (1- u) VI-\gamma _1V- \mu V]
+\lambda _3[\beta_0 (1- u) SI +\varepsilon \beta_0 (1- u) VI-\gamma I - \mu I]
\end{array}
\end{equation}
then, imposing first order conditions to minimize the Hamiltonian $H$ 
\begin{equation}
 \displaystyle{\frac{\partial H}{\partial u}}=2bu +I[\beta_0 S(\lambda _1-\lambda _3)+
 \varepsilon \beta_0 V(\lambda _2-\lambda _3)]=0,
  \end{equation}
we derive the optimal  restriction policy $u_{Q}^*$ (\ref{optcontrQdr}).\qed

\bigskip

Analogously to the quadratic case above, we can solve  the exponential case
\begin{corollary}
Let $\beta (u) = \beta_0 (1- u)$ and  $c_{exp}(u)=e^{k u}-1$.  

If $\lambda _3(t)>\lambda _1(t)$ and $\lambda _3(t)>\lambda _2(t)$, then the optimal control strategy $u_{exp}^*(t)$ for  problem (\ref{sistemar})--(\ref{funzionale}) is given by 
\begin{equation}
\label{optcontrexp}
u_{exp}^*(t)=
\left\{
\begin{array}{ll}
min\left \{max \left [0,\frac{1}{k}\ln \frac{\beta_0 I(t)K(t)}{k}\right ], \overline{u} \right \}
&
\mbox{if }K(t)>0
\\\\
0 & \mbox{ if }
K(t)\leq 0
\end{array}
\right.
\end{equation}
where K(t) is defined as
$\displaystyle{
K(t)=S(t)(\lambda _3(t)-\lambda _1(t))+
\varepsilon V(t)(\lambda _3(t)-\lambda _2(t))}
$.
\end{corollary}

\begin{remark}
We notice that in both cases, the optimal control strategy $u^*$ depends on the shadow price differences between infected and susceptible, and infected and vaccinated (see the paper Kruse \& Strack \cite{KS}). In other words, $\lambda _3-\lambda _1$ and $\lambda _3-\lambda _2$ can be interpreted as the marginal cost of having an additional susceptible person infected and as the marginal cost of having an additional vaccinated  person infected, respectively.

In the above Corollary, there is a sufficient condition so that the optimal control strategy $u_{exp}^*(\cdot)$ is well defined if the shadow price of the infected is higher than the shadow prices of the susceptible and vaccinee.
\end{remark}

\begin{remark}
It is interesting to point out that if in the optimal strategies $u_{Q}^*$ and $u_{exp}^*$ the maximum is not vanishing, then the optimal controls converges to the constant policy $\overline u$ as $b$ and $k$ tend to $0$ respectively. Hence if the ruler can cut off the social cost, then the optimal policy that can be adopted is the most strict one, represented by $\overline u$.
\end{remark}

\bigskip

Finally, we now consider the linear case: the proof of the following Proposition is reported in Appendix 1 and follows the reasoning in \cite{JL}. 

\begin{proposition} \label{optcontrlin_cor}
Let $\beta (u) = \beta_0 (1- u)$ and  $c_{lin}(u)=au$.  Then the optimal control strategy $u_{lin}^*(t)$ for  problem (\ref{sistemar})--(\ref{funzionale}) is given by
\begin{equation}
\label{optcontrlin}
u_{lin}^*(t)=
\left\{
\begin{array}{ll}
1
&
\mbox{if }\frac{\partial H}{\partial u}<0
\\\\
u_{sing}=\frac{A_2(t)}{A_1(t)}
&
\mbox{if }\frac{\partial H}{\partial u}=0
\\\\
0 & \mbox{if }\frac{\partial H}{\partial u}>0
\end{array}
\right.
\end{equation}
where $A_1 (\cdot)$ and $A_2 (\cdot)$ are defined in the proof.

\end{proposition}

\section{Numerical results} \label{numerics}

In this Section we apply the optimality results obtained in the previous section to analyze the controlled SVIR dynamic in several instances. In all our simulation results, we numerically solve the control problem by means of the forward-backward sweep algorithm (see \cite{LW}). This is a widely used indirect method for approximating the solution of optimal control problems (see e.g. \cite{MMH} for a convergence result). 

The parameters used in this section are summarized in Table (\ref{tabParam}), and they represent typical values of the recent Covid--19 pandemic. In particular, we fixed $\gamma^{-1}= 10.5$ and $\gamma_1^{-1}=14$ as the average of recovery time (from 7 to 14 days) and the time to reach full protection ($\approx 14$ days), respectively (see e.g. \texttt{https://www.who.int/health-topics/coronavirus}). The other parameters, $\beta_0$ and $\alpha$, has been estimated by using the aggregated Italian data provided by the ``Dipartimento della Protezione Civile'', available on GitHub\footnote{\texttt{https://github.com/pcm-dpc/COVID-19} and \texttt{https://github.com/italia/covid19-opendata-vaccini}.}, by exploiting the procedure described in the Appendix 2. Specifically, the transmission rate $\beta_0$ was obtained by using data from February, $24^{th}$ to March, $15^{th}$ 2020, setting $\alpha=\gamma_1=0$ (no vaccination and no containment measures at that time).  On the other hand, $\alpha$ was estimated with data from January $5^{th}$ (extension of third dose administration to people aged $\geq 12$), up to June, $8^{th}$, 2022. 

\smallskip

As introduced in Section \ref{sect_control}, the cost functional (\ref{funzionale}) is given by the sum of three terms, each related to a specific aspect of the problem: the ``social cost'' $J_{SC}(u) = \int_0^T c(u(t)) dt$, the ``infection cost'' $J_{IC}(u) = \int_0^T c_1 I(t) dt$, and the ``vaccination cost'' $J_{VC}(u) = \int_0^T c_2 \alpha S(t) dt$. In order to quantify the relative weights of each term, we rely on a quantification of the cost related to the hospitalized patients as available in the paper by Marcellusi et al. \cite{MFS22} (Supplementary material), assuming that the average cost for vaccination is $15$\texteuro \ per person. In particular, we considered the average daily cost weighted by the total number of patients  hospitalized  with and without ICU and only ICU\footnote{These costs per hospitalized patients have been estimated considering a sample of 996 Covid-19 hospitalisations recorded in Policlinico Tor Vergata Hospital between 2nd March 2020 and 27th December 2020.}.  Finally we normalized the corresponding weights in such a way $c_1=1$, implying $c_2=0.02$. 

\begin{table}[h] 
\centering
\begin{tabular}{lll}
\hline
Parameter & Value & Description \\ \hline \hline
$\beta_0$     & 0.220  &   transmission rate\\
$\gamma$     & 0.095 & recovery rate from infected \\ 
$\gamma_1$ & 0.071 &  recovery rate from vaccinated\\ 
$\alpha$        & 0.004 &   vaccination rate \\
\hline
\end{tabular}
\caption{Basic values of the parameters for the uncontrolled system. } 
\label{tabParam}
\end{table}

In the first set of experiments, we considered a time horizon of $240$ days with initial conditions $S_0=0.85$, $I_0=0.15$ and $V_0=0.0$. The parameters of the dynamical system (\ref{sistemar}) are shown in Table \ref{tabParam}. Furthermore, in order to fix the value for the scenario parameter $\varepsilon$, we choose to estimate it by using the vaccine effectiveness\footnote{The vaccine effectiveness is defined as the percentage reduction in risk of disease among vaccinated persons relative to unvaccinated persons.} $VE$ for the booster dose, averaged on the three available vaccines, ChAdOx1 nCoV-19 (Astra-Zeneca), BNT162b2 (Pfitzer-BionTech) and mRNA-1273 (Moderna), as reported in \cite{A2022} (Table 3): this procedure implied the value $\varepsilon \equiv (1-VE) =  0.078$. Finally, the birth-death rate has been set equal to zero\footnote{The annual birth and death rate in 2021 for the Italian population was estimated as $6.7\times 10^{-3}$ and $12 \times 10^{-3}$, respectively (source: ISTAT).}. Notice that in this case, $R_0^C = 0$.

\begin{figure}[H]
\hspace{-1.2cm}
\includegraphics[width=16cm, height=10cm]{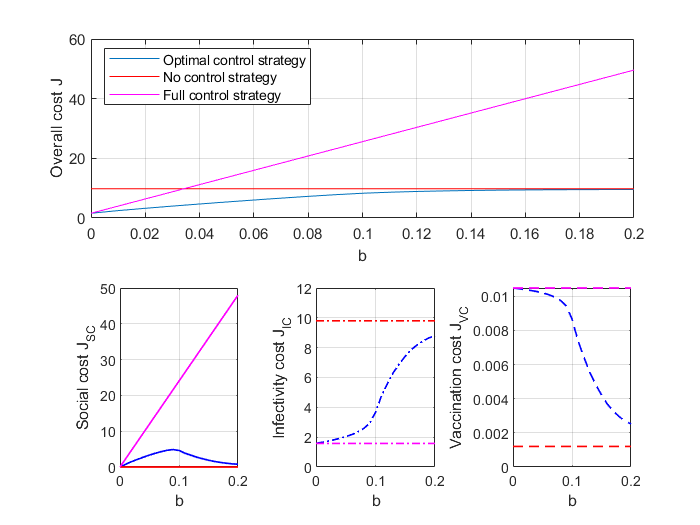}
\caption{Comparison of the cost values for the three strategies, no-control, full-control, optimal control, as a function of the social cost function parameter $b$.}
\label{Fig_costs_caso1}
\end{figure}

\paragraph{The quadratic cost function $c(u) = b u^2$.}
Firstly, we considered the quadratic cost function $c(u) = b u^2$ and studied the effect of the parameter $b$ on the cost functional, see Figure (\ref{Fig_costs_caso1}).  Increasing values of $b$ correspond to a higher social cost for increasing controls. As expected, we immediately see that if the social cost becomes larger ($b \nearrow$), the optimal strategy collapses to the no-control strategy: it would be convenient not to adopt any restrictions if the cost of such restrictions became too expensive. 

The ratio of the cumulative costs of each term in $J(u)$ with respect to the total cost is reported in Table (\ref{tabCosts1}). In the case of absence of control ($u(t) \equiv 0$), almost all of the expenditure is due to the disease, while for the completely controlled system ($u(t) \equiv 1$), the social cost amounts to slightly more than half of the total expenditure. In the optimally controlled system, this cost reduces to about 38\%.

\begin{figure}[H]
\hspace{-1.0cm}
\includegraphics[width=14cm, height=6cm]{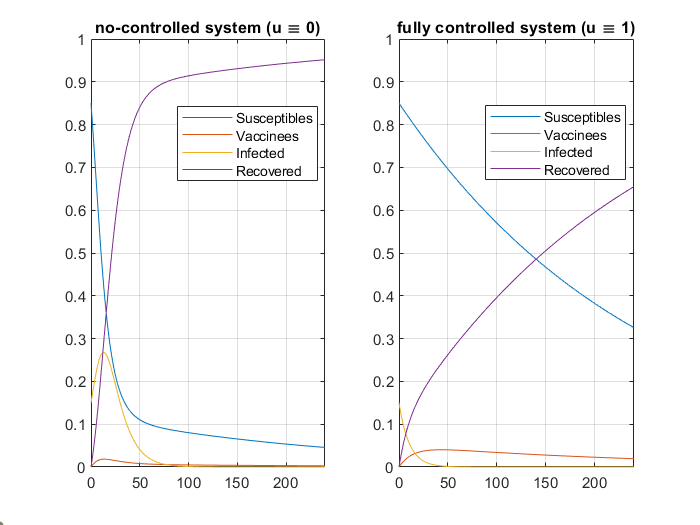}
\caption{Controlled and uncontrolled systems.}
\label{Fig_control1}
\end{figure}

\begin{figure}[H]
\hspace{-1.0cm}
\includegraphics[width=14cm, height=6cm]{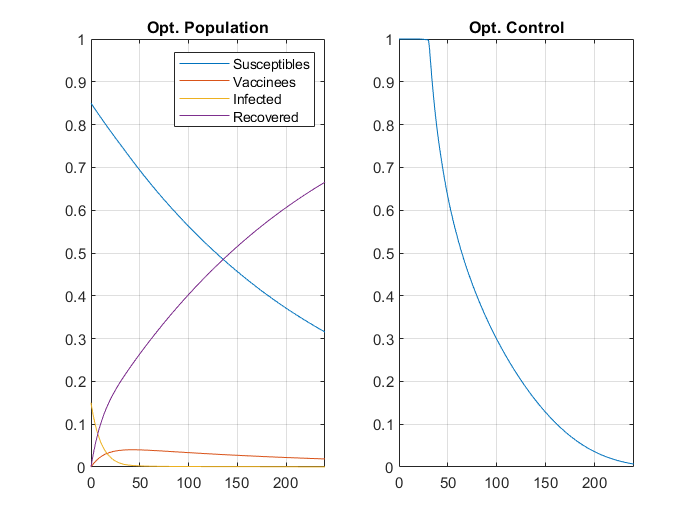}
\caption{Optimal controlled systems (on the left) and the optimal control strategy (on the right): quadratic cost function.}
\label{Fig_Optcontrol1}
\end{figure}

In Figure (\ref{Fig_control1}) we show the evolution of the state variables when no control is applied ($u(t) \equiv 0$ on the left) and with a constant full control ($u(t) \equiv 1$ on the right). The optimally controlled system is shown in Figure (\ref{Fig_Optcontrol1}), together with the optimal control function $u^*(t)$. We may notice, as expected,  that when the number of infected individuals is high, the optimal strategy is the highest level of control: hence, as this number is low enough, the optimal control switches to the lowest level. In this way the overall cost is reduced (see Table (\ref{tabCosts1})).

\begin{table}[h] 
\centering
\begin{tabular}{lc||ccc}
\hline
Control strategy & Total cost $J(u)$ &Social cost & Infection cost & Vaccination cost \\ \hline \hline
$u\equiv 0$ & 8.9264   &    0\%  & 99.9743\%   &  0.0257\% \\ 
$u\equiv 1$ & 6.3906   & 53.7734\% &   17.7020\% &   0.1175\% \\
$u^*        $ & 2.8705  & 38.4582\% &   61.1809\% &    0.3609\% \\
\hline
\end{tabular}
\caption{Total costs of the strategies and the relative social, infection and vaccination cost for the quadratic cost function. } 
\label{tabCosts1}
\end{table}

\begin{figure}[H]
\hspace{-1.5cm}
\includegraphics[width=15cm, height=9cm]{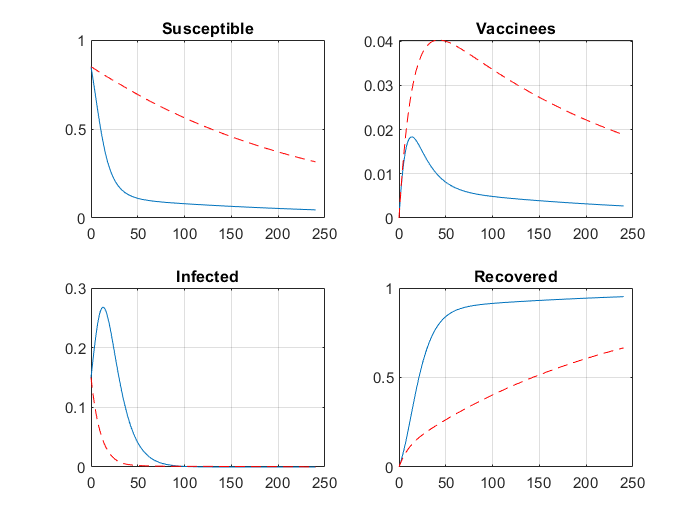}
\caption{Simulation results for the compartments dynamic with the social cost function $c(u) = 0.02 u^2$. Dashed lines are for the optimal control case, and solid lines for the "no-control" case.}
\label{Fig_OptCompart1}
\end{figure}

\begin{figure}[H]
\hspace{-1.2cm}
\includegraphics[width=16cm, height=10cm]{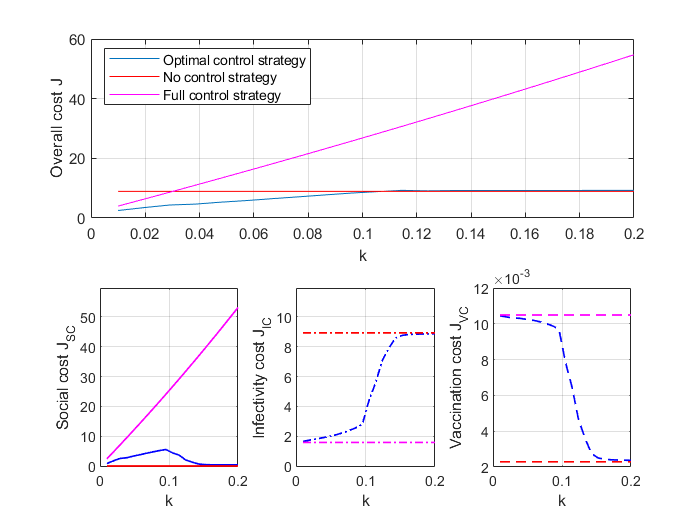}
\caption{Comparison of the cost values for the three strategies, no-control, full-control, optimal control, as a function of the exponential model of social cost function parameter $k$.}
\label{Fig_costs_caso2}
\end{figure}

\paragraph{The exponential cost function $c(u) = \exp(k u)-1$.}
In Figure (\ref{Fig_costs_caso2}) we plotted the values of the costs $J_{SC}, J_{IC}$ and $J_{VC}$ as a function of the parameter $k$. Similarly to the previous example, the social cost growth pushes the optimal strategy toward the no-control strategy. By choosing the value $k=0.06$, the optimal control strategy switches from the maximum control ($u \approx 1$) to the minimum one ($u \approx 0$) more quickly than in the quadratic case, but the effect on the optimal trajectories of the compartments is analogous. On the other hand, the reduction of the cost with respect to the full-control strategy is more pronounced. 

\begin{table}[h] 
\centering
\begin{tabular}{lc||ccc}
\hline
Control strategy & Total cost $J(u)$ &Social cost & Infection cost & Vaccination cost \\ \hline \hline
$u\equiv 0$ & 8.9264   &    0\%  & 99.9743\%   &  0.0257\% \\ 
$u\equiv 1$ & 16.4303 &  90.3257\%  &  9.6104\%   &  0.0639\% \\
$u^*        $ & 5.9897   & 64.2817\%   & 35.5483\% &   0.1700\% \\
\hline
\end{tabular}
\caption{Total costs of the strategies and the relative social, infection and vaccination cost for the exponential cost function. } 
\label{tabCosts2}
\end{table}

\begin{figure}[H]
\hspace{-1.0cm}
\includegraphics[width=14cm, height=6cm]{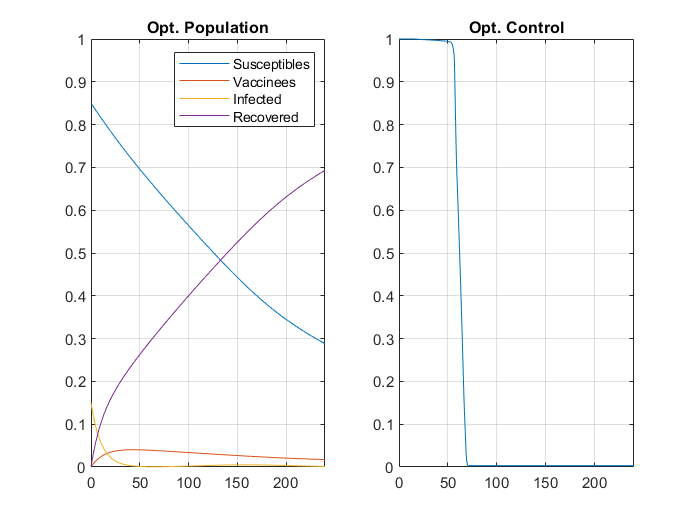}
\caption{Optimal controlled systems (on the left) and the optimal control strategy (on the right): exponential cost function.}
\label{Fig_Optcontrol2}
\end{figure}

\begin{figure}[H]
\hspace{-1.5cm}
\includegraphics[width=15cm, height=9cm]{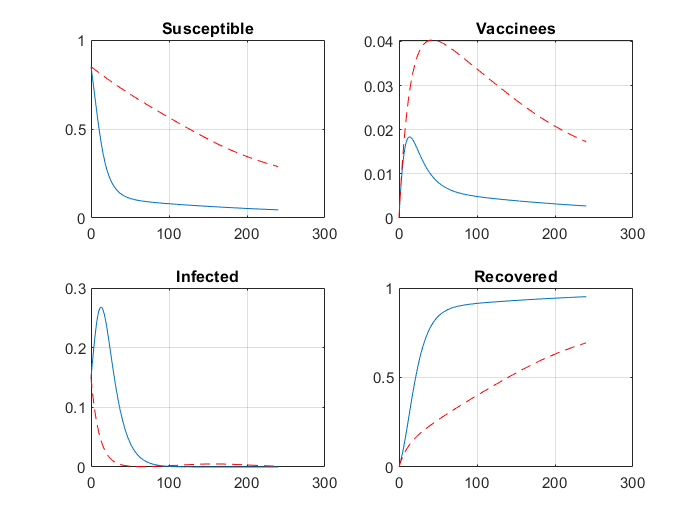}
\caption{Simulation results for the compartments dynamic with the social cost function $c(u) = \exp(0.06 u)-1$. Dashed lines are for the optimal control case, and solid lines for the "no-control" case.}
\label{Fig_OptCompart2}
\end{figure}

\begin{figure}[H]
\hspace{-1.2cm}
\includegraphics[width=16cm, height=10cm]{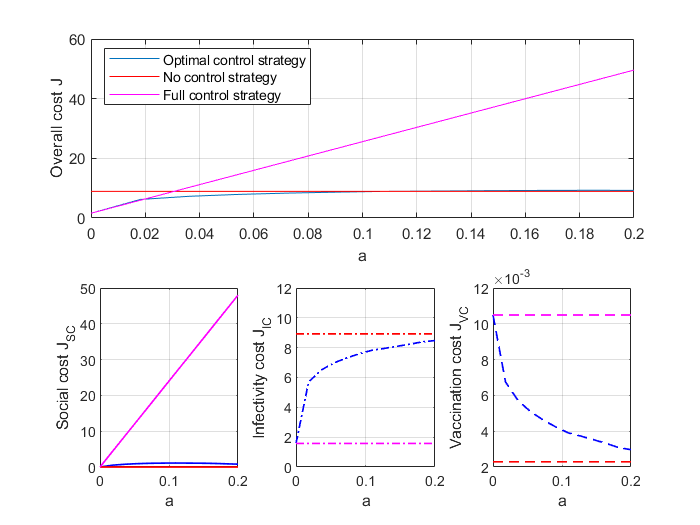}
\caption{Comparison of the cost values for the three strategies, no-control, full-control, optimal control, as a function of the linear model of social cost function parameter $a$.}
\label{Fig_costs_caso3}
\end{figure}

\paragraph{The linear cost function $c(u) =a u$.}
The dependence of the cost functional with respect to the parameter $a$ of the linear cost function is qualitatively similar to the other cases. An increase of the social cost pushes the optimal control towards the no-control strategy (see Figure (\ref{Fig_costs_caso3})). On the other hand, for a fixed value of $a$, the optimal control is ``bang-bang'', as proved in Proposition (\ref{optcontrlin_cor}), resulting in a different behavior of the optimal dynamic. As a matter of fact, the optimal trajectory of the infected compartment is decreasing in a first period when the optimal control is maximum, then the optimal control switches to the lowest value allowing the disease to spread again and then die out anyway. We observe, however, that the number of the infected nevertheless remains lower overall than the no-control case. Let us finally notice that the optimal, but discontinuous control generates a ``wave'' behavior of the Infected compartment.  

\begin{table}[h] 
\centering
\begin{tabular}{lc||ccc}
\hline
Control strategy & Total cost $J(u)$ &Social cost & Infection cost & Vaccination cost \\ \hline \hline
$u\equiv 0$ & 8.9264   &    0\%  & 99.9743\%   &  0.0257\% \\ 
$u\equiv 1$ & 13.5895  & 88.3034\%  & 11.6194\%   & 0.0772\% \\
$u^*        $ & 7.5057 & 14.6435\%  & 85.2786\%  &  0.0780\% \\
\hline
\end{tabular}
\caption{Total costs of the strategies and the relative social, infection and vaccination cost for the linear cost function. } 
\label{tabCosts3}
\end{table}

\begin{figure}[H]
\hspace{-1.0cm}
\includegraphics[width=14cm, height=6cm]{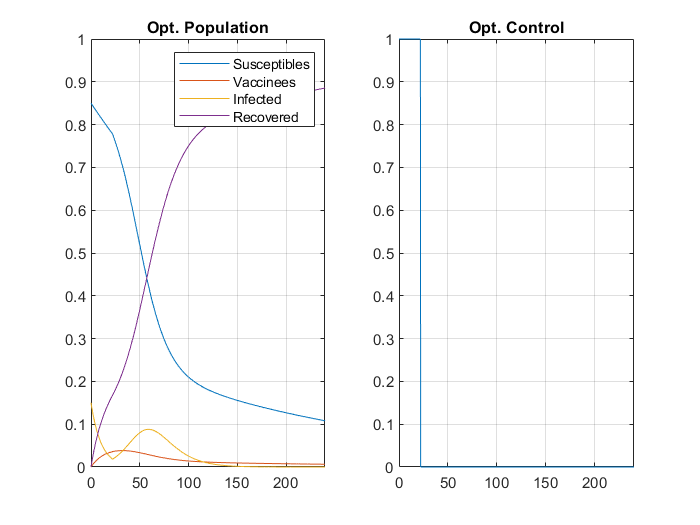}
\caption{Optimal controlled systems (on the left) and the optimal control strategy (on the right): linear cost function.}
\label{Fig_Optcontrol3}
\end{figure}

\begin{figure}[H]
\hspace{-1.5cm}
\includegraphics[width=15cm, height=9cm]{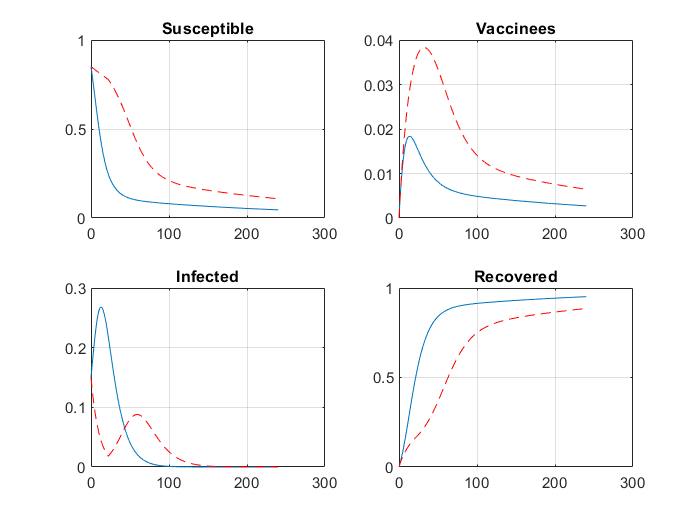}
\caption{Simulation results for the compartments dynamic with the social cost function $c(u) = 0.05u$. Dashed lines are for the optimal control case, and solid lines for the "no-control" case.}
\label{Fig_OptCompart3}
\end{figure}

\begin{figure}[H]
\hspace{-1.5cm}
\includegraphics[width=16cm, height=10cm]{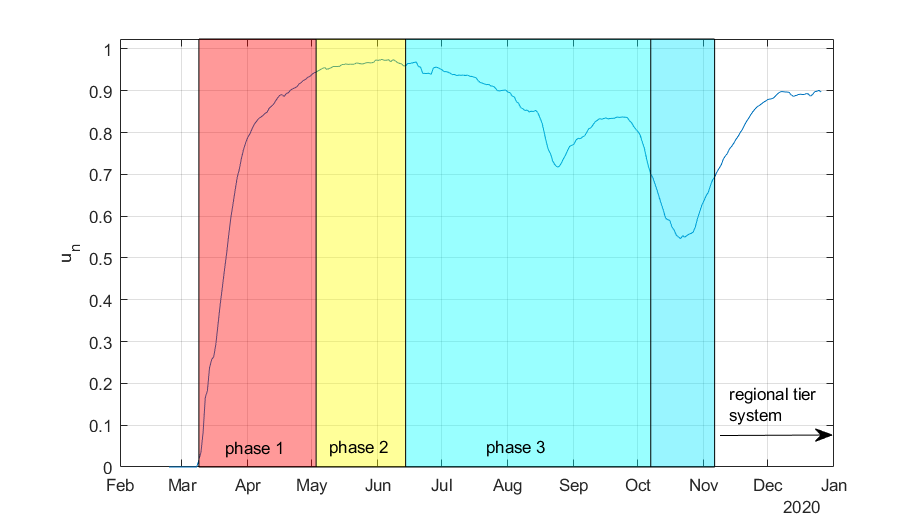}
\caption{Ex-post valuation of the policy $\hat u_n$ in the period February 24th - December 26th, 2020.}
\label{Fig_ex_post_val}
\end{figure}

\paragraph{Ex-post valuation of the realized strategy.}
As a further experiment, we report the ex-post valuation of the strategy adopted in Italy during the period February 24th - December 26th, 2020. Since the vaccination campaign started on December 27th, 2020, we considered the dynamical model with $\alpha=\gamma_1\equiv0$. By taking the discrete-time setting described in Appendix 2 for the model dynamic, and exploiting the linear relationship between the rate and the control variable, see (\ref{lin_contr}), we immediately get the estimate of the realized control variable.
$$
\hat u_n = 1 - \frac{\hat \beta_n}{\hat \beta_0},
$$
where $\hat \beta_0$ has been estimated as described at the beginning of the Section \ref{numerics}, and $\hat \beta_n$ is the day-by-day solution of the linear system built from the discrete-time SIR model (see Appendix 2). We estimated these parameters using the data from the GitHub repository, and the results are shown in Figure (\ref{Fig_ex_post_val}), where we highlighted the different phases of the pandemic event.

\begin{figure}[H]
\hspace{-1.0cm}
\includegraphics[width=14cm, height=6cm]{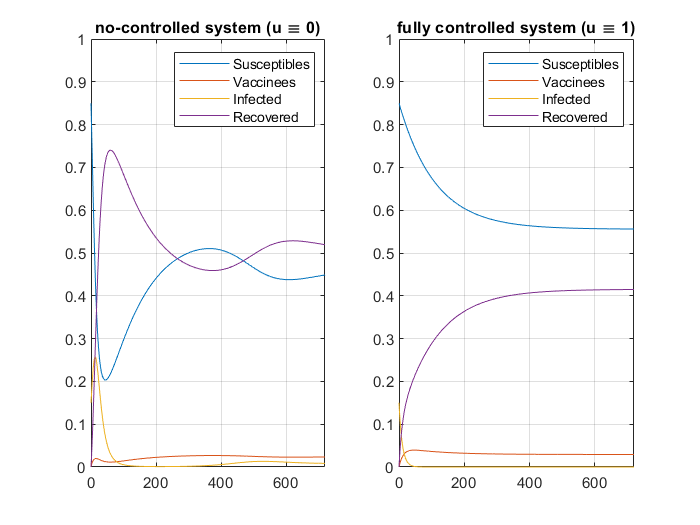}
\caption{The uncontrolled and fully controlled systems with the quadratic social cost function. Here $R_0^C=1.22$.}
\label{Fig_endemic01}
\end{figure}

\begin{figure}[H]
\hspace{-1.0cm}
\includegraphics[width=14cm, height=6cm]{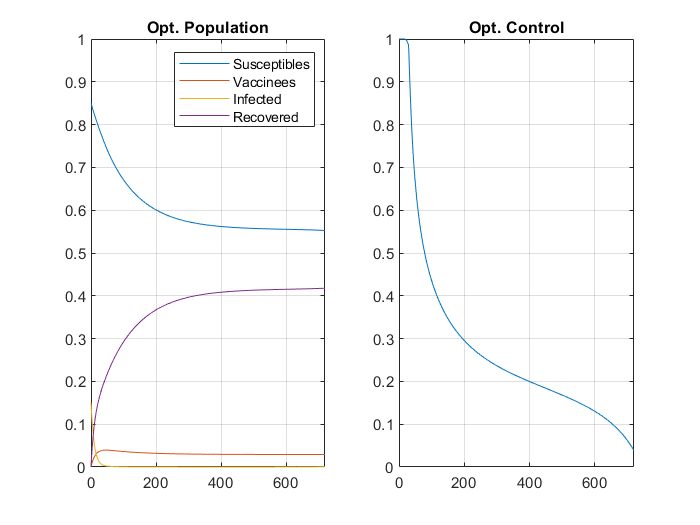}
\caption{Optimal controlled systems (on the left) and the optimal control strategy (on the right): quadratic cost function.}
\label{Fig_endemic_control}
\end{figure}

\begin{figure}[H]
\hspace{-1.0cm}
\includegraphics[width=15cm, height=9cm]{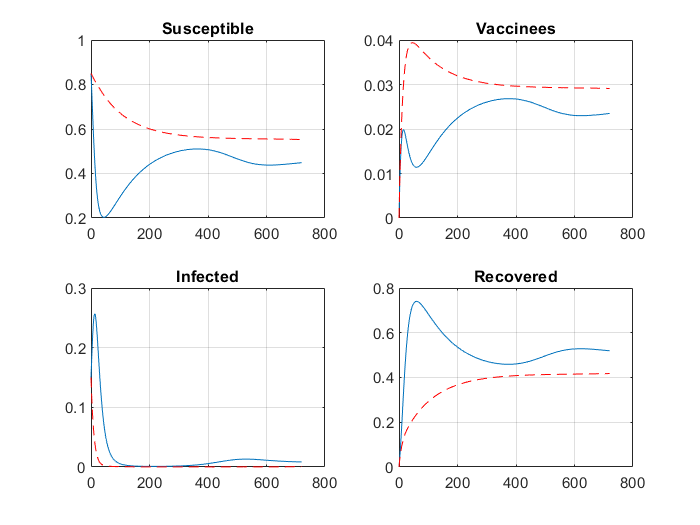}
\caption{Simulation results for the compartments dynamic with the quadratic social cost function. Dashed lines are for the optimal control case, and solid lines for the "no-control" case.}
\label{Fig_endemic_compar}
\end{figure}

\paragraph{An example with $R_0^C>1$.} As a final set of experiments, we report the analysis of the control problem for a SVIR model showing an endemic equilibrium (see Section \ref{basic_svir}). With the same set of parameters as before (see Table \ref{tabParam}), it is sufficient to set the (daily) birth-death rate equal to $0.005$ to get $R_0^C=1.22$. The uncontrolled and fully controlled systems are shown in Figure (\ref{Fig_endemic01}), in the case of a quadratic social cost function. It is immediately seen, that in this example $I_+ > 0$. The solution of the optimal control problem is shown in Figures (\ref{Fig_endemic_control}) and (\ref{Fig_endemic_compar}) and the Table (\ref{tabCosts5}) reports the corresponding costs. We may immediately see the effectiveness of the optimal control strategy in this example. 


\begin{table}[h] 
\centering
\begin{tabular}{lc||ccc}
\hline
Control strategy & Total cost $J(u)$ &Social cost & Infection cost & Vaccination cost \\ \hline \hline
$u\equiv 0$ & 13.3261   &    0\%  & 99.8121\%    & 0.1879\% \\ 
$u\equiv 1$ & 15.9359  & 90.3622\%  & 9.4206\%   & 0.2172\% \\
$u^*        $ & 3.4631 & 50.5445\%  & 48.4597\%  & 0.9957\% \\
\hline
\end{tabular}
\caption{Total costs of the strategies and the relative social, infection and vaccination cost for the quadratic cost function in the endemic example.} 
\label{tabCosts5}
\end{table}

\section{Concluding remarks}

In this paper we considered the problem of optimal control of an infectious disease by modeling its diffusion through a SVIR compartmental dynamic. Differently from the classical SIR (or SIRD) model, we consider the possibility to implement a vaccination campaign to immunize the population. The control of the disease is realized by adopting political measures of containment, generally identified here as \textit{social distancing}, since they may have an impact on the social behavior of the population, e.g., the use of face masks, the partial or total closure of many activities, educational structures, commercial activities, production, and/or different degrees of prohibition of movements.
All these measures aim to reduce the disease's diffusion but imply a cost for the whole society. Hence we introduced a cost functional that explicitly considers these social measures, the social cost function, other than the cost of vaccination and the cost due to the infected population.
Our main result is the characterization of the solution of the optimal control of the SVIR dynamical model obtained in terms of a controlled diffusion rate to minimize the overall cost of the disease. We thoroughly describe the optimal control strategy using the Pontryagin Maximum Principle for several instances of the social cost function. Finally, we implemented the optimal controlled system using the Forward-Backward Sweep algorithm by calibrating some of the model parameters using the Italian dataset of the recent Covid--19 pandemic.  

As a final contribution, we highlight some possible further developments of this research. First of all, it would be natural to include in the control variables the rate of vaccination, other than the control of the rate of diffusion. The addition of a controlled vaccination campaign makes the problem more complex from a mathematical point of view, resulting in a two-dimensional constrained optimization problem. In this paper we preferred to focus only on the impact of the social measures over the overall cost. Nevertheless, we deserve to include this variable in future research. As a second remark, in light of the recent pandemic, it would be interesting to introduce some modification of the basic compartmental dynamical model. In particular, for the SVIR model we considered, the possibility for a recovered to be re-infected. This could be relevant to capturing the phenomenon of virus mutations. The resulting dynamical model can be simply obtained by adding a link from the $R$ to the $I$ compartment. Of course, the analysis of the stability of the uncontrolled system is quite different from the one presented here, as well as the corresponding Hamiltonian function. 

Finally, other essential SVIR dynamical system modifications are still under consideration to meet the peculiarities of a possible pandemic. In particular, as noticed in the numerical section based on the Covid--19 data, the cost of the pandemic can be particularly severe for its impact on hospital services becoming overwhelmed, and it is "quantifiable". On the other hand, the cost of being "quarantined" is undoubtedly more challenging to assess. It is, therefore, sensible to split the "I" compartment into at least three sub-compartments, e.g., Hospitalized w/o ICU, Hospitalized with ICU, and Quarantined, and to adjust the cost function accordingly, eventually adding other compartments as the Death and/or the Exposed, thus resulting in the SVIRD/SVEIR/SVEIRD models.

%
%
%
%
%

\bigskip

\textbf{Acknowldgments:} The Authors would like to thank Prof. F.S. Mennini for having provided some of the empirical data, and Prof. R. Cerqueti for the many helpful suggestions. All errors remain ours.
%

\section*{Appendix 1} \textbf{Proof of Corollary} (\ref{optcontrlin_cor}). The system (\ref{sistemar}) can be written as
\begin{equation}
\label{sistemarl}
\left\{
\begin{array}{ll}
S' =\beta _0(u-1)SI -(\alpha + \mu )S+\mu &S(0)=S_0
\\\\
V' =\alpha S+\varepsilon \beta _0(u-1)VI-(\gamma _1+\mu )V &V(0)=V_0
\\\\
I' =\beta _0(1-u)SI+\varepsilon \beta _0(1-u)VI- (\gamma +\mu )I &I(0)=I_0
\end{array}
\right.
\end{equation}
The Hamiltonian function $H$ is defined as
\begin{equation}
\label{hamiltonianal}
\begin{array}{l}
H(t,S,V,I,u) 
=au+c_1I+c_2\alpha S+
\lambda _1[-\beta_0(1- u) S I -\alpha S+ \mu (1- S) ]+
\\\\
+\lambda _2[\alpha S+\varepsilon \beta_0 (u-1) VI-(\gamma _1+\mu )V]
\\\\
+\lambda _3[\beta_0 (1- u) SI +\varepsilon \beta_0 (1- u) VI-(\gamma I + \mu )I]
\end{array}
\end{equation}
and the co--state system (\ref{costato}) can be written as
\begin{equation}
\label{costatol}
\left\{
\begin{array}{l}
\lambda _1' =\beta _0(1-u)(\lambda _1 - 
\lambda _3)I+\alpha (\lambda _1 - 
\lambda _2)+\mu \lambda _1-c_2\alpha 
\\\\
\lambda _2' =\varepsilon\beta _0(1-u)(\lambda _2 -\lambda _3)I+(\gamma _1+\mu )\lambda _2
\\\\
\lambda _3' =\beta _0(1-u)(\lambda _1 - 
\lambda _3)S+\varepsilon\beta _0(1-u)(\lambda _2 -\lambda _3)V+(\gamma +\mu )\lambda _3-c_1
\end{array}
\right.
\end{equation}
Since the Hamiltonian is linear in the control $u$ is bang--bang, singular or a combination, see \cite{JL}, \cite{LW}. The singular case is attained if 
\begin{equation}
 \displaystyle{\frac{\partial H}{\partial u}}=a +\beta_0 (\lambda _1-\lambda _3)SI+
 \varepsilon \beta_0 (\lambda _2-\lambda _3)VI=0,
  \end{equation}
  on a non--trivial interval of time, else if $\frac{\partial H}{\partial u}<0$ the optimal control would be at its upper bound, while if $\frac{\partial H}{\partial u}>0$ it would be ta its lower bound.
  So, to study the singular case, suppose $\frac{\partial H}{\partial u}=0$ on a non--trivial interval of time, calculate
  $$
  \frac {d}{dt}\left ( \frac{\partial H}{\partial u}\right )=0
  $$
  and show that, in the above equation, the control $u$ does no appear explicitly. Then, the value of the singular control will be obtain evaluating

$$
  \frac {d^2}{dt^2}\left ( \frac{\partial H}{\partial u}\right )=0
  $$
Hence
$$
\begin{array}{l}
\displaystyle{
 0= \frac {d}{dt}\left ( \frac{\partial H}{\partial u}\right )=}
 \\\\
 =(\lambda _1'-\lambda _3')SI+(\lambda _1-\lambda _3)S'I+(\lambda _1-\lambda _3)SI'+
 \\\\
 +\varepsilon (\lambda _2'-\lambda _3')VI+\varepsilon (\lambda _2-\lambda _3)V' I+\varepsilon (\lambda _2-\lambda _3)VI'
 \end{array}
  $$
We calculate each term of the sum separately and then we add them together:
$$
\begin{array}{l}
(\lambda _1'-\lambda _3')SI=
[\beta _0(1-u)(\lambda _1 - 
\lambda _3)I+\alpha (\lambda _1 - 
\lambda _2)+\mu \lambda _1-c_2\alpha -
\\\\
-\beta _0(1-u)(\lambda _1 - 
\lambda _3)S-\varepsilon\beta _0(1-u)(\lambda _2 -\lambda _3)V-(\gamma +\mu )\lambda _3+c_1
]SI
\end{array}
$$
\vspace{0.5cm}
$$
(\lambda _1-\lambda _3)S'I=
[\beta _0(u-1)SI -(\alpha + \mu )S+\mu
](\lambda _1-\lambda _3)I
$$
\vspace{0.5cm}
$$
(\lambda _1-\lambda _3)SI'=
[\beta _0(1-u)SI+\varepsilon \beta _0(1-u)VI- (\gamma +\mu )I](\lambda _1-\lambda _3)S
$$
\vspace{0.5cm}
$$
\begin{array}{l}
\varepsilon (\lambda _2'-\lambda _3')VI=
\varepsilon [\varepsilon\beta _0(1-u)(\lambda _2 -\lambda _3)I+(\gamma _1+\mu )\lambda _2-
\\\\
-\beta _0(1-u)(\lambda _1 - 
\lambda _3)S-\varepsilon\beta _0(1-u)(\lambda _2 -\lambda _3)V-(\gamma +\mu )\lambda _3+c_1
]VI
\end{array}
$$
\vspace{0.5cm}
$$
\varepsilon (\lambda _2-\lambda _3)V' I=
\varepsilon [\alpha S+\varepsilon \beta _0(u-1)VI-(\gamma _1+\mu )V](\lambda _2-\lambda _3) I
$$
\vspace{0.5cm}
$$
\varepsilon (\lambda _2-\lambda _3)VI'=
\varepsilon [\beta _0(1-u)SI+\varepsilon \beta _0(1-u)VI- (\gamma +\mu )I](\lambda _2-\lambda _3)V
$$
Hence, summing up we obtain
$$
\frac {d}{dt}\left ( \frac{\partial H}{\partial u}\right )=
$$
$$
=\alpha (\lambda _1 - 
\lambda _2)SI+\mu \lambda _1SI-c_2\alpha SI-(\gamma +\mu )\lambda _3SI+c_1SI
-(\alpha + \mu )(\lambda _1-\lambda _3)SI+\mu(\lambda _1-\lambda _3)I-
$$
$$
- (\gamma +\mu )(\lambda _1-\lambda _3)SI
+\varepsilon (\gamma _1+\mu )\lambda _2VI
-\varepsilon (\gamma +\mu )\lambda _3VI+\varepsilon c_1VI-
$$
$$
-\varepsilon (\gamma _1+\mu )(\lambda _2-\lambda _3)VI-\varepsilon (\gamma +\mu )(\lambda _2-\lambda _3)VI=
$$
$$
=[-(\gamma +\mu)\lambda _1+(\alpha +\mu )\lambda _3-\alpha \lambda _2+(c_1-\alpha c_2)]SI+\mu(\lambda _1-\lambda _3)I+
$$
$$
+\varepsilon [-(\gamma +\mu)\lambda _2+(\gamma _1+\mu )\lambda _3+c_1]VI
$$

Since the control doesn't appear in the previous expression, we compute the second derivative\footnote{Computations were realized with the help of Mathematica\copyright}:
$$
\frac {d^2}{dt^2}\left ( \frac{\partial H}{\partial u}\right )=
$$
$$
= I' (S (c_1-\alpha c_2-\alpha  \lambda_2+\alpha  \lambda_3-(\gamma +\mu ) \lambda_1+\mu \lambda_3)+\varepsilon  V
(c_1-(\gamma +\mu ) \lambda_2+(\gamma_1+\mu ) \lambda_3)+\mu (\lambda_1-\lambda_3))+
$$
$$
+ I (S' (c_1-\alpha c_2-\alpha \lambda_2+\alpha \lambda_3-(\gamma +\mu ) \lambda_1+\mu \lambda_3)+\varepsilon  V'(c_1-(\gamma +\mu ) \lambda_2+(\gamma_1 +\mu) \lambda_3)+
$$
$$
+ S  (-\alpha \lambda_2'+(\alpha +\mu ) \lambda_3' - ((\gamma +\mu ) \lambda_1' ) )+\varepsilon  V ((\gamma_1+\mu ) \lambda_3'-(\gamma +\mu ) \lambda_2' )+\mu  (\lambda_1'- \lambda_3')).
$$
By substituting (\ref{sistemar}) and (\ref{costato}), we get
$$
\frac {d^2}{dt^2}\left ( \frac{\partial H}{\partial u}\right )=
$$
$$
= I ((S (\alpha +\beta  I-\beta  I u+\mu )-\mu ) (\alpha  \lambda_2-\alpha \lambda_3+\gamma  \lambda_1-c_1+\alpha c_2+\lambda_1 \mu-\lambda_3 \mu ) +
$$
$$
S ((\alpha+\mu ) (\gamma \lambda_3 - c_1+\lambda_3 \mu-\beta  \lambda_1 S(u-1)+\beta \lambda_3 S (u-1) - \beta  \varepsilon \lambda_2 (u-1) V+\beta \varepsilon  \lambda_3 (u-1) V)-
$$
$$
(\gamma +\mu ) (-\alpha \lambda_2+\alpha (-c_2)+\lambda_1(\alpha +I (\beta -\beta  u)+\mu)+\beta  I \lambda_3(u-1))-\alpha  \lambda_2(\gamma_1-\beta  \varepsilon I (u-1)+\mu )+
$$
$$
\alpha  \beta \varepsilon  (-I) \lambda_3(u-1))+\mu  (-\alpha  \lambda_2+c_1-\alpha c_2+\lambda_1 (\alpha+I (\beta -\beta  u)+\mu )+\beta  I\lambda_3 (u-1)+\beta \lambda_1 S(u-1)-
$$
$$
\lambda_3 (\gamma +\mu+\beta  S (u-1)+\beta  \varepsilon  u V-\beta  \varepsilon  V)+\varepsilon\beta  \lambda_2 (u-1)V)-(\gamma +\mu +\beta  S (u-1)+\beta \varepsilon  u V-\beta  \varepsilon V) \times
$$
$$
(S (-\alpha  \lambda_2+\alpha  \lambda_3-\lambda_1 (\gamma +\mu)+c_1-\alpha c_2+\lambda_3 \mu)+\varepsilon  V (-\lambda_2 (\gamma +\mu )+\lambda_3 (\gamma_1+\mu)+c_1)+\mu  (\lambda_1-\lambda_3))+
$$
$$
\varepsilon (-\lambda_2 (\gamma +\mu)+\lambda_3 (\gamma_1+\mu )+c_1) (\alpha  S-V(\gamma_1-\beta  \varepsilon I (u-1)+\mu ))+\varepsilon  V((\gamma_1+\mu ) (\gamma \lambda_3-c_1+\lambda_3 \mu-
$$
$$
\beta  \lambda_1 S(u-1)+\beta  \lambda_3 S(u-1)-\beta  \varepsilon \lambda_2 (u-1) V+\beta \varepsilon  \lambda_3 (u-1)V)-(\gamma +\mu ) (\gamma_1\lambda_2-\beta  \varepsilon I (u-1) (\lambda_2-\lambda_3)+\lambda_2 \mu ))),
$$
which is linear in the control $u$: hence
$$
\frac {d^2}{dt^2}\left ( \frac{\partial H}{\partial u}\right )=A_1(t) u(t) - A_2(t) = 0
$$
giving the singular control
$$
u_{sing}(t) = \frac{A_2(t)}{A_1(t)}
$$
if $A_1(t) \neq 0$ and $0 \leq \frac{A_2(t)}{A_1(t)} \leq \bar u$, where the functions $A_i$ are given by
$$
A_1 = I (-\alpha \lambda_1 \mu +2 \alpha \lambda_2 \mu -\alpha \lambda_3 \mu +2 \gamma \lambda_1 \mu + c_1 (S (-\alpha (\varepsilon -2)+\gamma +\beta I+3 \mu -2 \beta \varepsilon V)+
$$
$$
\varepsilon V (\gamma +2 \gamma_1+\beta \varepsilon I+3 \mu )-2 \mu -\beta S^2-\beta \varepsilon ^2 V^2)+\alpha  c_2 (-S (\alpha +2 \gamma +\beta I+3 \mu -\beta \varepsilon V)+2 \mu +\beta S^2)-
$$
$$
\beta I \lambda_1 \mu +\beta I \lambda_3 \mu +\alpha \beta \varepsilon I \lambda_2 S-\alpha \beta \varepsilon I \lambda_3 S-\alpha \beta I \lambda_2 S+\alpha \beta I \lambda_3 S-\beta \gamma I \lambda_3 S-
$$
$$
\beta \gamma \varepsilon ^2 I \lambda_3 V+\beta \gamma_1 \varepsilon ^2 I \lambda_3 V+ \lambda_1 \mu ^2- \lambda_3 \mu ^2-\alpha \beta \lambda_1 S^2+\alpha \beta \lambda_2 S^2+\beta \gamma \lambda_1 S^2-\alpha ^2 \lambda_2 S+\alpha ^2 \lambda_3 S+
$$
$$
\alpha \gamma \varepsilon \lambda_2 S-2 \alpha \gamma \lambda_2 S-\alpha \gamma_1 \varepsilon \lambda_3 S+\alpha \gamma_1 \lambda_2 S+\alpha \varepsilon \lambda_2 \mu S-\alpha \varepsilon \lambda_3 \mu S-
$$
$$
2 \alpha \lambda_2 \mu S+2 \alpha \lambda_3 \mu S-\gamma ^2 \lambda_1 S-2 \gamma \lambda_1 \mu S- \lambda_1 \mu ^2 S+ \lambda_3 \mu ^2 S+\beta \gamma \varepsilon \lambda_1 S V+\beta \gamma \varepsilon \lambda_2 S V-
$$
$$
\beta \gamma_1 \varepsilon \lambda_1 S V+\beta \gamma \varepsilon ^2 \lambda_2 V^2-\beta \gamma_1 \varepsilon ^2 \lambda_2 V^2-\beta \varepsilon \lambda_1 \mu V-\gamma ^2 \varepsilon \lambda_2 V-
$$
$$
2 \gamma \varepsilon \lambda_2 \mu V+\gamma_1^2 \varepsilon \lambda_3 V+2 \lambda_1 \varepsilon \lambda_3 \mu V-\varepsilon \lambda_2 \mu ^2 V+\varepsilon \lambda_3 \mu ^2 V+\varepsilon \beta \lambda_2 \mu V),
$$
and
$$
A_2 = I (\beta c_1 (I (S+\varepsilon ^2 V )-(S+\varepsilon V)^2 )+\alpha \beta c_2 S (-I+S+\varepsilon V)-\beta I \lambda_1 \mu +\beta I \lambda_3 \mu +\alpha \beta \varepsilon I \lambda_2 S-
$$
$$
\alpha \beta \varepsilon I \lambda_3 S-\alpha \beta I \lambda_2 S+\alpha \beta I \lambda_3 S-\beta \gamma I \lambda_3 S-\beta \gamma \varepsilon ^2 I \lambda_3 V+\beta \gamma_1 \varepsilon ^2 I \lambda_3 V-
$$
$$
\alpha \beta \lambda_1 S^2+\alpha \beta \lambda_2 S^2+\beta \gamma \lambda_1 S^2+\beta \gamma \varepsilon \lambda_1 S V+\beta \gamma \varepsilon \lambda_2 S V-\beta \gamma_1 \varepsilon \lambda_1 S V+
$$
$$
\beta \gamma \varepsilon ^2 \lambda_2 V^2-\beta \gamma_1 \varepsilon ^2 \lambda_2 V^2-\beta \varepsilon \lambda_1 \mu V+\varepsilon \beta \lambda_2 \mu V ).
$$

\qed

\section*{Appendix 2}

The estimation of parameters for the SVIR model is typically based on a discrete-time version of (\ref{svirconst}) (where we set $\beta_1 = \epsilon \beta$), once available the observations for the compartments. By considering a discretization period $\Delta t =1$ day, and by letting $n = 0, 1, \ldots$ the discrete time instants, we consider a first-order scheme to obtain the following finite-differences model:
\begin{equation} \label{discrSvir}
\left\{
\begin{array}{ll}
S_{n+1} =S_n -\beta S_n I_n -\alpha S_n+ \mu - \mu S_n & S_0=s_0
\\\\
V_{n+1} =V_n + \alpha S_n - \epsilon \beta  V_n I_n-\gamma _1V_n - \mu V_n & V_0 = v_0
\\\\
I_{n+1} =I_n + \beta S_n I_n + \epsilon \beta V_n I_n-\gamma I_n - \mu I_n & I_0=i_0
\\\\
R_{n+1}=R_n + \gamma _1 V_n+\gamma I_n- \mu R_n & R_0=r_0.
\end{array}
\right.
\end{equation}
By assuming the parameter $\mu$ as exogeneously given, the system is linear in the remaining  parameters $\vartheta := (\beta, \alpha, \gamma_1, \gamma)'$. Hence, by defining 
$$
\Delta_n = \left( \begin{array}{c}
S_{n+1}-S_n(1-\mu) \\
V_{n+1}-V_n(1-\mu) \\
I_{n+1} - I_n(1-\mu) \\
R_{n+1}-R_n(1-\mu)
\end{array}
\right) \ \ \ A_n = \left(\begin{array}{cccc} 
- S_n I_n & -S_n & 0 & 0 \\
\epsilon V_n I_n & S_n & -V_n & 0 \\
(S_n+\epsilon V_n) I_n & 0 & 0 & -I_n \\
0 & 0 & V_n & I_n
\end{array} \right)
$$
we can write (\ref{discrSvir}) in matrix form as  $A_n \vartheta = \Delta_n$. Given the observed values of the compartments, the standard constrained regression OLS can therefore be used as the basic estimation procedure for the parameters of the model, that is
$$
\hat \vartheta = \mbox{argmin}_{\vartheta \geq 0} \sum_{n=1}^{T-1} ||\Delta_n - A_n \vartheta ||_2^2.
$$
(see e.g. \cite{CNP20}). 

By considering the system instead with time-varying coefficients
\begin{equation} \label{discrSvir}
\left\{
\begin{array}{ll}
S_{n+1} =S_n -\beta_n S_n I_n -\alpha_n S_n+ \mu - \mu S_n & S_0=s_0
\\\\
V_{n+1} =V_n + \alpha_n S_n - \epsilon \beta_n  V_n I_n-\gamma_{1,n}V_n - \mu V_n & V_0 = v_0
\\\\
I_{n+1} =I_n + \beta_n S_n I_n + \epsilon \beta_n V_n I_n-\gamma_n I_n - \mu I_n & I_0=i_0
\\\\
R_{n+1}=R_n + \gamma_{1,n} V_n+\gamma_n I_n- \mu R_n & R_0=r_0,
\end{array}
\right.
\end{equation}
we may write the discrete-time dynamic equations in a matrix form as $A_n \vartheta_n = \Delta_n$,
where $\vartheta_n := (\beta_n, \alpha_n, \gamma_{1,n}, \gamma_n)'$. When $\alpha_n = \gamma_{1,n} \equiv 0$, the system has the unique solution
$$
\left\{ \begin{array}{ccc}
\hat \beta_n & = & - \frac{\Delta_n^{(1)}}{S_n I_n} \\ \\
\hat \gamma_n & = & \frac{\Delta_n^{(4)}}{I_n}.
\end{array} \right.
$$


\begin{thebibliography}{99}
\bibitem{A} Abakuks A. An Optimal Isolation Policy for an Epidemic. \textit{Journal of Applied Probability} 10(2), 247--262, (1973)

\bibitem{AAL21} Alvarez F. D. Argente and F. Lippi. A Simple Planning Problem for COVID-19 Lock-down, Testing, and Tracing. \textit{American Economic Review: Insights} 3(3), 367--82, (2021)

\bibitem{A2022} N. Andrews et al. Covid-19 Vaccine Effectiveness against the Omicron (B.1.1.529) Variant. \textit{The new England Journal of Medicine}, 386, (2022)

\bibitem{Ben00} Behncke, H. Optimal control of deterministic epidemics. \textit{Optim. Control Appl. Methods} 21 (6), 269--285, (2000)

\bibitem{B} Bolzoni L.,  Bonacini E., Soresina C., Groppi M.
Time-optimal control strategies in SIR epidemic models, \textit{Mathematical Biosciences} 292, 86--96, (2017)

\bibitem{BC} Brauer, F. and Castillo-Chavez, C. Mathematical Models in Population Biology and Epidemiology. Springer Science, Berlin (2010) 

\bibitem{CNP20} Calafiore G.C., Novara C., Possieri C. A time-varying SIRD model for the COVID-19 contagion in Italy. \textit{Annual Reviews in Control} 50, 361--372,  (2020)

\bibitem{CGLZ22} Calvia, A., Gozzi, F., Lippi, F., Zanco, G. A Simple Planning Problem for COVID-19 Lockdown: A Dynamic Programming Approach. Preprint available online at \texttt{https://arxiv.org/pdf/2206.00613.pdf}  (2022)

\bibitem{FF} Federico S., Ferrari G. Taming the spread of an epidemic by lockdown policies, \textit{Journal of Mathematical Economics}, 93, 102453, (2021)

\bibitem{FFT22} Federico S., Ferrari G., Torrente M.L. Optimal vaccination in a SIRS epidemic model, Preprint available online at \texttt{https://arxiv.org/abs/2206.03284}, (2022)

\bibitem{FR} Fleming W.H. and Rishel R.W. Deterministic and stochastic optimal control, Applications of Mathematics, vol. 1, Springer-Verlag, New York, Heidelberg and Berlin, (1975)

\bibitem{GS} Gaff, H., \& Schaefer, E. Optimal control applied to vaccination and treatment strategies for various epidemiological models. \textit{Mathematical Biosciences \&  Engineering}, 6(3), 469--492, (2009).

\bibitem{HW73} Hethcote, H.W., Waltman, P. Optimal Vaccination Schedules in a Deterministic Epidemic Model. \textit{Mathematical Biosciences} 18, 365--381, (1973).

\bibitem{JL} Joshi, H. R., Lenhart, S., Hota, S., \& Augusto, F. B. Optimal control of an SIR model with changing behavior through an education campaign. \textit{Electronic Journal of Differential Equations} 50, 1--14, (2015)

\bibitem{KM} Kermack, W.O., McKendrick, A.G. Contributions to the mathematical theory of epidemics. \textit{Bulletin of Mathematical Biology} 53 (1--2), 33--55, (1991)

\bibitem{KS} Kruse, T., Strack, P. Optimal control of an epidemic through social distancing. Preprint at SSRN: https://ssrn.com/abstract=3581295, (2020)

\bibitem{KPS} Kumar, A., \& Srivastava, P. K. Vaccination and treatment as control interventions in an infectious disease model with their cost optimization. \textit{Communications in Nonlinear Science and Numerical Simulation} 44, 334---343,  (2017).

\bibitem{LS}  Ledzewicz U., Schättler H. On optimal singular controls for a general SIR-model with vaccination and treatment. Conference Publications, 2011 (Special) , 981--990, (2011).

\bibitem{LW} Lenhart, S. and Workman, J.T. Optimal Control Applied to Biological Models. Mathematical and Computational Biology Series, Chapman \& Hall/CRC, London, (2007)

\bibitem{XTI} Liu, Xianning, Yasuhiro Takeuchi, and Shingo Iwami. SVIR epidemic models with vaccination strategies. \textit{Journal of Theoretical Biology} 253(1), 1--11, (2008)

\bibitem{MFS22} Marcellusi A., Fabiano G., Sciattella P., Andreoni M., Mennini F.S. The Impact of COVID--19 Vaccination on the Italian Healthcare System: A Scenario Analysis. \textit{Clinical Drug Investigation} 42, 237--242, (2022)

\bibitem{MMH} McAsey M, Mou L, Han W. Convergence of the forward-backward sweep method in optimal control. \textit{Computational Optimization and Applications} 53, 207--226, (2012)

\bibitem{MSW20} Miclo, L., Spiroz, D., Weibull, J. Optimal Epidemic Suppression under an ICU Constraint. Preprint available online at \texttt{https://arxiv.org/abs/2005.01327}, (2020)

\bibitem{WMV} Witbooi, Peter J., Grant E. Muller, and Garth J. Van Schalkwyk.  Vaccination control in a stochastic SVIR epidemic model. \textit{Computational and mathematical methods in medicine}, doi: 10.1155/2015/271654 (2015)

\end{thebibliography}
\end{document}